\documentclass[11pt]{article}
\usepackage{graphicx}
\usepackage{amssymb}
\usepackage{epstopdf}
\usepackage{epsfig}
\usepackage{subfigure}
\usepackage{amsmath, amssymb, amsfonts}
\usepackage{a4, a4wide}
\usepackage{lscape}
\usepackage{color}
\usepackage{verbatim}
\usepackage{array}
\usepackage{hyperref}
\usepackage{arydshln}

\def\ie{{\rm i.e.,\/}\ }

\newcommand{\nco}{\newcommand}
\nco{\one}{\ensuremath{\,\,\mathrm{l}\!\!\!1}} 
\nco{\NN}{\mathbb{N}}
\nco{\ZZ}{\mathbb{Z}}
\nco{\QQ}{\mathbb{Q}}
\nco{\RR}{\mathbb{R}}
\nco{\CC}{\mathbb{C}}
\nco{\HH}{\mathbb{H}}
\nco{\OO}{\mathbb{O}}
\nco{\red}{\color{red}}
\nco{\redend}{\normalcolor}
\nco{\colorend}{\normalcolor}
\nco{\0}{\underline{0}}
\nco{\1}{\underline{1}}
\nco{\2}{\underline{2}}
\nco{\3}{\underline{3}}
\def\ommit#1{{}}
\makeatletter
\newcommand{\pmatrixcmd}[1]{}%
\DeclareRobustCommand{\pmatrixcmd}[1]{\left(\env@matrix#1\endmatrix\right)}
\makeatother
\begin{document}

\noindent
 \vspace{0.2cm} 
\begin{center}
\begin{Large}
    {Character tables (modular data) for Drinfeld doubles of finite groups}
\vspace{0.5cm}
\end{Large}
\renewcommand{\thefootnote}{\arabic{footnote}}
{Seminar presented at the 7th International Conference on Mathematical Methods in Physics\\ Joint Conference: CBPF-IMPA-ICTP-SISSA-TWAS, Rio de Janeiro, April 2012} 
\end{center}
\date{\vspace{-5ex}}
\begin{center}
\renewcommand{\thefootnote}{\arabic{footnote}}
{R. Coquereaux} \footnotemark[1]${}^,$ \footnotemark[2]
\\ \renewcommand{\thefootnote}{\arabic{footnote}} 
\end{center}


\abstract{In view of applications to conformal field theory or to other branches of theoretical physics and mathematics, 
new examples of character tables for Drinfeld doubles of finite groups (modular data) are made available on a website.}

\vspace{0.3cm}



\addtocounter{footnote}{0}
\footnotetext[1]{{\scriptsize{\it IMPA, Instituto Nacional de Matem\'atica Pura e Aplicada, Rio de Janeiro, Brasil}.}}
\footnotetext[2]{{\scriptsize{\it UMI CNRS-IMPA 2924, on leave from CPT, UMR 7332, Luminy, Marseille, France}.}}

\vspace{0.3cm}

\unitlength = 1mm

\section{Introduction}

\subsection{Summary}
  
  To any finite group $H$ one can associate not only its group algebra $\CC H$ (a cocommutative Hopf algebra) but also a non-cocommutative Hopf algebra $D(H)$ called its Drinfeld double \cite{Drinfeld:quasibialgebras}. 

  Usual representation theory on $\CC$ of a finite group is the representation theory of its group algebra. One can also consider the representation theory of its Drinfeld double.
  Fusion rules describing the reduction of a tensor product of irreducible representations of  $\CC H$ can be deduced from the character table of the group.
    In the same way,  fusion rules describing the reduction of a tensor product of irreducible representations of $D(H)$ can be deduced from the modular $S$ matrix of the double, which is often called the ``finite group modular data'', like in \cite{CGR:modulardata}, or, sometimes, the ``character table of its quantum double''.
      
  One can associate with $D(H)$ a fusion category (of representations) that is modular -- \ie we have an action of the modular group $SL(2,\ZZ)$ generated by the so-called $S$ and $T$ modular matrices on the linear span of the irreps.
  {The construction of $D(H)$ can be twisted by a cocycle, but here we only refer to the untwisted case.}

  Character tables for finite groups can be found in many places. This is not so for their modular data.
  Although general formulae for $S$ and $T$ exist in the literature \cite{CGR:modulardata}, \cite{DiPaRo:double}, and although the subject is not really new, 
  very few explicit results -- \ie very few tables of modular data-- are available, even in the untwisted situation (see the discussion in section \ref{sec:availableresults}).
  In view of applications to conformal field theory,  and for possible applications of Drinfeld doubles to the theory of quantum computers  \cite{E.Rowell:quantumcomputing},  \cite{MThompson:quantumcomputing}, or maybe for other reasons  yet to be discovered, such tables should be calculated and made accessible. Because of their size, $S$ matrices cannot usually be published in printed form, even for  small groups.
The main  purpose of this contribution is to announce that several new examples are now available on-line,  on the web site: \\  \url{http://www.cpt.univ-mrs.fr/~coque/quantumdoubles/comments.html}

Our motivation, for undertaking these calculations, was to provide the necessary material that would allow one to test several conjectures generalizing the results obtained in 
 \cite{RCJBZ:sumrules}. \\ We refer the interested reader to the forthcoming article \cite{RC-JBZ:doublesWIP}.

 For the finite groups listed on the website -- a list that includes all exceptional finite subroups of $SU(2)$ and $SU(3)$, and a few members from their infinite series of subgroups--  the available tables give their modular data, \ie the matrices $S$ and $T$, 
 and in most cases the  fusion matrices; each example also comes with a summary text file. 
Information about the conjugacy classes and their centralizers was  obtained from GAP \cite{GAP}. 
Part of the calculation was performed with Magma \cite{Magma}, with an algorithm which, like in \cite{MThompson:quantumcomputing},Ê uses a variant of equation \ref{eq: GannonSmatrix}.
The results were then handled to Mathematica \cite{Mathematica}. On the website, 
$S$ matrices are given as lists of lists whereas  $T$ matrices and fusion matrices are usually described as sparse arrays (Mathematica syntax). 
Calculation of fusion matrices, using the Verlinde formula \cite{Verlinde}, was performed with Mathematica. 
Properties of the Drinfeld doubles of  finite subgroups of $SU(2)$ and $SU(3)$, together with their fusion graphs, are given and commented, together with a number of other topics,  in the article \cite{RC-JBZ:doublesWIP}.

The present contribution contains: \\
1) A sketch of the Drinfeld double construction -- see \cite{DiPaRo:double, Drinfeld:quasibialgebras, KoornEtAl} -- that we repeat here mostly for cultural reasons,  since we are only interested in the representation ring of the latter (in our framework, knowing ``what'' is represented is not necessarily useful).\\
 2) The formulae giving $S$ and $T$ which can be extracted from \cite{CGR:modulardata} or \cite{DiPaRo:double}, but we use a variant of the latter (eq. 3).\\
 3) Details about one example whose modular data can be found on the aforementioned website, namely the Drinfeld double of the binary tetrahedral group.
 
 \subsection{Miscellaneous remarks}
 \label{sec:availableresults}
 
  \paragraph{About existing tables.}      
      As it was mentioned, character tables for finite groups can be found in many books and computerized data basis.  They are used by many scientific communities, ranging from chemists and crystallographers to mathematicians.  
    In the same way, for applications of quantum doubles (in particular Drinfeld doubles) one would need to have access to their modular data\footnote{We already mentioned that the $S$ matrix is the analog of a character table, but the latter is usually normalized in such a way that the first line is $1,1....1$, so that a better analog of the character table is not $S$ itself but the matrix $\chi$ with elements $\chi(m,n) = S_{n,m} / S_{n,1}$.} (in particular to the $S$ matrix) without having to calculate it from scratch.
    Unfortunately, such tables for Drinfeld doubles of finite groups, abstractly defined about $25$ years ago, are not easily available.
      The situation is even worse for more general quantum doubles (ie twisted Drinfeld doubles)  and for $c\neq 0$ orbifold theories where one considers a Lie group $G$ at level $k$ together with a finite subgroup $H$ of $G$.
     Of course, general expressions that we shall recall later giving the $S$ and $T$ modular generators for Drinfeld doubles can be found in \cite{CGR:modulardata} and  \cite{DiPaRo:double}, but these are intricate formulae involving summations over characters of various subgroups. These are not ``tables''.
      Several explicit examples of $S$ matrices for $c=1$ orbifolds have been worked out for instance in \cite{CappDAppo} and \cite{KacTodorov}, and several explicit formulae for $S$ matrices of Drinfeld doubles can be found in \cite{CGR:modulardata}, in the case of abelian groups (in particular  cyclic groups $\ZZ_n$),  dihedral groups $D_n$, and for the permutation group $S_3$.  More recently, and in the framework of investigations on quantum computing,  fusion rules were published by \cite{MThompson:quantumcomputing} for the doubles of $D_3$, $S_3$, $S_4$ and of the alternating group $A_4$.
            
      \paragraph{About motivations.}   
      Several results (sum rules) recently obtained by \cite{RCJBZ:sumrules} for WZW theories, more precisely for representation theory of simple Lie groups at level $k$ (Kac-Moody algebras) and at infinite level, \ie for the Lie groups themselves, lead to various conjectures that we wanted to test on the representation theory of Drinfeld doubles -- see the article \cite{RC-JBZ:doublesWIP}. The lack of available explicit results  for modular data concerning doubles made necessary their determination for a variety of cases.  Such results cannot be published in printed form;  indeed, for instance, the rank of the Drinfeld double of the subgroup $\Sigma_{360 \times 3}$ of $SU(3)$, a group with class number $17$,  is $240$, so that its $S$ matrix is a  $240 \times 240$  matrix with complex coefficients. 
      
\section{A sketch of the Drinfeld double construction}

Le $H$ be a finite group, $F(H)$, the (commutative) algebra of functions on $H$.
A base is  $\{\delta_g\}$ where $g$ are  the group elements.
Let $C(H)$ the (cocommutative) group algebra of $H$.
A base is $\{g\}$,  better call it  $\{x\}$.
The Drinfeld double, as a vector space, is : $D(H) = F(H) \otimes C(H)$.
A base is $\{\delta_g  \otimes x\}$.  
Multiplication is $(\delta_g \otimes x)( \delta_h \otimes y) = \delta_{g,xhx^{-1}} \, \delta_g  \otimes xy$ .
Comultiplication  is $\Delta ( \delta_g \otimes x ) =  \sum_{{h,k \in H ,  hk=g}} (\delta_h \otimes x) (\delta_k \otimes x)$.
Counit is $\epsilon( \delta_g \otimes x ) = \delta_{g,e} \, 1 \otimes e$.
Antipode is ${\mathcal S}( \delta_g \otimes x) = \delta_{x^{-1}g^{-1}x} \otimes x^{-1}$.
The Drinfeld double $D(H)$ is a quasi-triangular Hopf algebra (existence of an R matrix obeying Yang-Baxter equations).
Its R-matrix is
$R = \sum_{g\in H} (\delta_g \otimes e) \bigotimes (1\otimes g) \; \in D(H) \bigotimes D(H)$.
More generally one can construct the Drinfeld double of a Hopf algebra: the result is a quasi-triangular Hopf algebra, even if the first is not.
The square of the antipode (an automorphism) is given by conjugation by an invertible element $u$:
$ \forall a \in D(H), \, {\mathcal S}^2(a) = u a u^{-1}$.
The category of representations of $D(H)$ is modular
 (it was not so for the category of representations of $C(H)$).
There is a finite number of simple objects (irreducible representations or irreps).
Call $i,j$ two simple objects.
The modular group $SL(2,Z)$, generated by two elements called  $S$ and $T$, acts: 
one build matrices $S_{ij}$ and $T_{ij}$ indexed by irreps $i,j$ of $D(H)$ and obeying
 $C=S^2=(ST)^3, C^2=S^4=1$.
 To an irrep $j$ of  $D(H)$ one associates a complex number (called conformal weight $h_j$).
 To this (finite list) of numbers $h_j$ we associate the (diagonal) $T$ matrix:
 $diag(T) =  exp[- 2 i \pi \, c/24]  \; exp[2 i \pi \,  \{h_i\}] $, here we take $c=0$.
 From $u$ (using the square of the antipode ${\mathcal S}$) and the list $\{h_j\}$ we associate the matrix $S$ :
 $ S_{ij}/S_{00} = e^{2 i \pi (h_i + h_j)} \; tr_{i \otimes j} \, \Delta(u)$.
Now consider tensor products of irreducible representations of $D(H)$. 
They are usually reducible:  $i \otimes j = \sum_k \, N_{ij}^k \;  k$.
From the $S$ matrix one can determine (Verlinde formula) the fusion coefficients $N_{ij}^k$ (non negative integers).
 To $i$, we associate a fusion matrix $N_i$ with elements $(N_i)_j^k= N_{ij}^k$.
 To the matrix $N_i$ we  associate a fusion graph whose adjacency matrix is $N_i$.

\section{General formulae for $S$ and $T$}
\label{eq: GannonSmatrix}

As discussed in \cite{CGR:modulardata, DiPaRo:double, Drinfeld:quasibialgebras}  there is a one to one correspondence between irreps of  the Drinfeld double $D(H)$ and pairs
$(A, \alpha)$ where
$A$ is a conjugacy class of the finite group $H$, 
and $\alpha$ is an irrep of the centralizer (defined up to conjugation), in $H$, of any
representative\footnote{Two elements belonging to the same conjugacy class have conjugated  (hence isomorphic) centralizers.} element of $A$.

  
Call $C_a$ the centralizer, in $H$, of a representative element $a$ of the conjugacy class $A$. Notice that $\vert A \vert \, \vert C_a \vert = \vert H \vert$.
The following expressions for $S$ and $T$ matrix elements of the untwisted case can be extracted from \cite{CGR:modulardata} or \cite{DiPaRo:double}.
Let $A$ and $B$ be two conjugacy classes of $H$.  We take $a \in A$, $b \in B$, call  $C_a$ and $C_b$ their centralizers,  and choose
$\alpha$, an irrep of $C_a$, and $\beta$, an irrep of $C_b$. Then, 
\begin{eqnarray}
S_{(A,\alpha)(B,\beta)} &=& \frac{1}{| C_a||C_b|} \sum_{\scriptstyle{ g\in H}\atop \scriptstyle{ a \, g b g^{-1}=g b g^{-1}a}}
\chi_{\alpha}(g b g^{-1})^* \chi_{\beta}(g^{-1}ag)^*\\
T_{(A,\alpha)(B,\beta)} &=& \delta_{AB}\, \delta_{\alpha\beta}\, \frac{\chi_{\alpha}( a)}{\chi_{\alpha}(e)}\,,
\end{eqnarray}
where $\chi_\alpha$ and $\chi_\beta$ are the irreducible characters associated with
the irreps $\alpha$ and $\beta$ of the groups $C_a$ and $C_b$. The neutral element is $e$.

Equivalently, let ${\mathcal T}_a=\{a_i\}$  (resp. ${\mathcal T}_b=\{b_j\}$) be a system of coset representatives for the {left}  classes of  $H/C_a$ (resp. a system of  coset representatives for the left classes of  $H/C_b$), then, 

\begin{eqnarray}
S_{(A,\alpha)(B,\beta)}&=& \frac{1}{\vert H \vert}  \sum_{g_{ij}}
\chi_{\alpha}(g_{ij} b g_{ij}^{-1})^* \chi_{\beta}(g_{ij}^{-1}ag_{ij})^*
\end{eqnarray}
where the sum runs over all 
$g_{ij} = a_i b_j^{-1}$, with $a_i \in {\mathcal T}_a$, $b_j \in {\mathcal T}_b$, that  obey $[b_j^{-1} b b_j,a_i^{-1} a a_i ] = 1$;  here $[ {}\,,\, {}]$ is  the commutator defined as $[a,b]=a^{-1}b^{-1}ab$.

\section{Double of the binary tetrahedral group (example)}

The tetrahedral group is defined in $SO(3)$ as the finite subgroup of order $12$ preserving a tetrahedron.
Obviously it is isomorphic with $A_4$, the alternating group on $4$ objects (even permutations of the $4$ vertices of a tetrahedron).
The binary tetrahedral group $H$, or order $24$, is its two-fold cover in $SU(2)$.
It is isomorphic with  $SL(2,F_3)$, the group of $2 \times 2$ matrices with entries in the field $F_3$, and determinant one.
Here is the list of the 24 elements:
{\tiny
\[
\begin{array}{c}
\left(
\begin{array}{cc}
 0 & 1 \\
 2 & 0 \\
\end{array}
\right),\left(
\begin{array}{cc}
 0 & 1 \\
 2 & 1 \\
\end{array}
\right),\left(
\begin{array}{cc}
 0 & 1 \\
 2 & 2 \\
\end{array}
\right),\left(
\begin{array}{cc}
 0 & 2 \\
 1 & 0 \\
\end{array}
\right),\left(
\begin{array}{cc}
 0 & 2 \\
 1 & 1 \\
\end{array}
\right),\left(
\begin{array}{cc}
 0 & 2 \\
 1 & 2 \\
\end{array}
\right),\left(
\begin{array}{cc}
 1 & 0 \\
 0 & 1 \\
\end{array}
\right),\left(
\begin{array}{cc}
 1 & 1 \\
 2 & 0 \\
\end{array}
\right),
\cr
\left(
\begin{array}{cc}
 1 & 2 \\
 1 & 0 \\
\end{array}
\right),\left(
\begin{array}{cc}
 2 & 0 \\
 0 & 2 \\
\end{array}
\right),\left(
\begin{array}{cc}
 2 & 1 \\
 2 & 0 \\
\end{array}
\right),\left(
\begin{array}{cc}
 2 & 2 \\
 1 & 0 \\
\end{array}
\right),\left(
\begin{array}{cc}
 1 & 0 \\
 1 & 1 \\
\end{array}
\right),\left(
\begin{array}{cc}
 1 & 0 \\
 2 & 1 \\
\end{array}
\right),\left(
\begin{array}{cc}
 1 & 1 \\
 0 & 1 \\
\end{array}
\right),\left(
\begin{array}{cc}
 1 & 1 \\
 1 & 2 \\
\end{array}
\right),
\cr
\left(
\begin{array}{cc}
 1 & 2 \\
 0 & 1 \\
\end{array}
\right),\left(
\begin{array}{cc}
 1 & 2 \\
 2 & 2 \\
\end{array}
\right),\left(
\begin{array}{cc}
 2 & 0 \\
 1 & 2 \\
\end{array}
\right),\left(
\begin{array}{cc}
 2 & 0 \\
 2 & 2 \\
\end{array}
\right),\left(
\begin{array}{cc}
 2 & 1 \\
 0 & 2 \\
\end{array}
\right),\left(
\begin{array}{cc}
 2 & 1 \\
 1 & 1 \\
\end{array}
\right),\left(
\begin{array}{cc}
 2 & 2 \\
 0 & 2 \\
\end{array}
\right),\left(
\begin{array}{cc}
 2 & 2 \\
 2 & 1 \\
\end{array}
\right)
\end{array}
\]
}

The character table of $H$ is given below:   
 {\footnotesize
\[
\begin{array}{c:c|ccccccc}
Classes &
& 
{\tiny \begin{pmatrix} 1 & 0 \\ 0 & 1 \end{pmatrix}} 
& 
{\tiny \begin{pmatrix} -1 & 0 \\ 0 & -1 \end{pmatrix}} 
&
{\tiny \begin{pmatrix} 0 & 1 \\ -1 & -1 \end{pmatrix}} 
&
{\tiny \begin{pmatrix} -1 & -1 \\ 1 & 0 \end{pmatrix}} 
&
{\tiny \begin{pmatrix} 0 & 1 \\ -1 & 0 \end{pmatrix}} 
&
{\tiny \begin{pmatrix} 0 & 1 \\ -1 & 1 \end{pmatrix}} 
&
{\tiny \begin{pmatrix} 1 & -1 \\ 1 & 0 \end{pmatrix}} 
\\
\hline
numbering & & 1 & 2 & 3 & 4 & 5 & 6 & 7 \\
\hline
\#(class) & &{  1} & 1 & 4 & 4 & 6 & {  4} & {  4}  \\
\hline
order & & {  1} & 2 & 3 & 3 & 4 & 6 & 6 \\
\hline
{\tiny irreps} \\  {\tiny numbering \, ,\, names} & & & & & & &\\
1, \quad\alpha_0& +& {  1}&1  &1 &1 &{  1} &{  1} &{  1} \\
\hdashline
2, \quad \alpha_1&0 &{  1}& 1&-1- j  &j &{  1} &{  j} &{  -1-j}  \\
3, \quad \alpha_5& 0 &{  1} &1  & j &-1-j &{  1} &{  -1-j} &{  j} \\
4, \quad\alpha_6&- &2 & -2 &-1 &-1 &0 &1 &1 \\
5, \quad\alpha_2&0 &2 &-2  &1+j &-j &0 &j &-1-j \\
6, \quad\alpha_4&0 &2 &-2  &-j &1+j & 0&-1-j &j \\
7, \quad\alpha_3&+ &{  3} & 3 &0 &0 &{  -1} &{  0}&{  0}  \\
\hline
C_g & & {  H}& H & Z_2 \times Z_3 & Z_2 \times Z_3 & {  Z_2} \times {  Z_2} & Z_2 \times {  Z_3} & Z_2 \times { Z_3} \\
\vert \widehat C_g \vert& & {  7} & 7 & 6 & 6 & 4 &  6 & 6  \\

\end{array}
\]
}
where 
$j = exp(2 i \pi/3)$ and the types of irreps are denoted as follows: 
$+$ real type,
$0$ complex type, 
$-$ quaternionic type.
The fourth line of the table lists the order $p$ of the elements (in the sense $g^p=1$) in each conjugacy class.
The penultimate line gives the centralizers of the conjugacy classes,  and the last line gives the class numbers of the centralizers (the number of their irreps). All other entries should be clear.
The irreps of $H$ are named according to the components of the highest root of $E_6$ on the base of coroots,
since the dimensions of the non trivial irreps of $H$ are the coroots labels of $E_6$:  $\{1, 2, 3, 2, 1, 2\}$, i.e.   $\theta =  1\alpha_1 + 2 \alpha_2  + 3 \alpha_3 + 2 \alpha_4 + 1\alpha_5 + 2 \alpha_6$. We call 
 $\alpha_0$ the trivial. The embedding of $H$ into $SU(2)$ is faithfully realized  by the two dimensional representation $\alpha_6$ (the fundamental representation).
 
 The character table for the tetrahedral group itself (not its binary), that we do not use here,  could be obtained from the above by restricting the table to lines (irreps) numbered 1,2,3,7 and to columns (conjugacy classes) 1,5,6,7;
 its classes have respectively $(1,3,4,4)$ elements, with orders $(1,2,3,3)$ and centralizers $A_4$, $\ZZ_2 \times \ZZ_2$, $\ZZ_3$ and $\ZZ_3$.
 
In order to calculate the $S$ and $T$ matrices of the quantum double, we  need the character tables of the centralizers, for  group representatives of all conjugacy classes. 
This is quite easy in the present case since centralizers of the two classes $1$ and $-1$ are given by the group $H$ itself, whereas the centralizers of the five others are abelian groups. 
$H$ possesses $7$ inequivalent irreps, and abelian groups have a number of irreps equal to their degree.

So we a priori know that the $S$ and $T$ matrices can be decomposed in $7$ blocks of sizes ${7,7,6,6,4,6,6}$, and that the total number of irreps of the Drinfeld double is  $42$.
Calculating the matrix elements of $S$ and $T$ is now a straightforward (but cumbersome) task. 
$S$ being symmetric, it is enough to display the blocks $S[I,J]$ for $J \geq I$. 
The results are given at the end of this section.

One recovers the (usual) character table of $H$ from its modular data by first restricting the matrix $S$ to the lines of its first block (first $7$ lines) and by selecting the first column of each consecutive block;
the expected result (up to a global conjugation) is obtained by further dividing each column (labelled by a conjugacy class) by the cardinal of its class $(1,1,4,4,6,4,4)$, and by multiplying everything by $\vert H \vert$.
Indeed,  $S_{(A=[e], \alpha), (B,1)} = 1/(\vert C_e \vert) \, 1/(\vert C_B \vert) \, \sum_g \chi_\alpha(g b g^{-1})^*$  simplifies to $(\vert B \vert /\vert H \vert) \,  \chi_\alpha(b)^*$ since $\chi$ is then central for $H$, 
and since $1/\vert C_B \vert = \vert B \vert /\vert H \vert$.

The quantum dimensions $\mu_j$ of the irreps $j$ defined as $S_{1,j}/S_{1,1}$, that we give for the consecutive blocks,
 are $\{1, 1, 1, 2, 2, 2, 3\}$,  $\{1, 1, 1, 2, 2, 2, 3\}$,  $\{4, 4, 4, 4, 4, 4\}$,   $\{4, 4, 4, 4, 4, 4\}$,  $\{6, 6, 6, 6\}$,  $\{4, 4, 4, 4, 4, 4\}$, $\{4, 4, 4, 4, 4, 4\}$,
and the global dimension of the Drinfeld double $\sum_j \, \mu_j^2$, equal to $1/\vert S_{1,1} \vert^2$ by unitarity, is $24^2$. 
The global dimension of the Drinfeld double of a finite group is actually always equal to the square of the order of the latter.

Fusion matrices $N_j$ can now be obtained from $S$ by using the Verlinde formula. These matrices have non-negative integer coefficients and can be considered as adjacency matrices of graphs (the fusion graphs).
We do not display the $42$ fusion graphs but only the one associated with the fundamental irrep $\alpha_6$, which is the fourth of the first block, as already mentioned.  The fourth of the second block gives a similar graph since
this block is again associated with $H$ itself, now viewed as the centralizer of $-1$.

\begin{figure}[htp]
  \begin{center}
 \includegraphics[width=8.0cm]{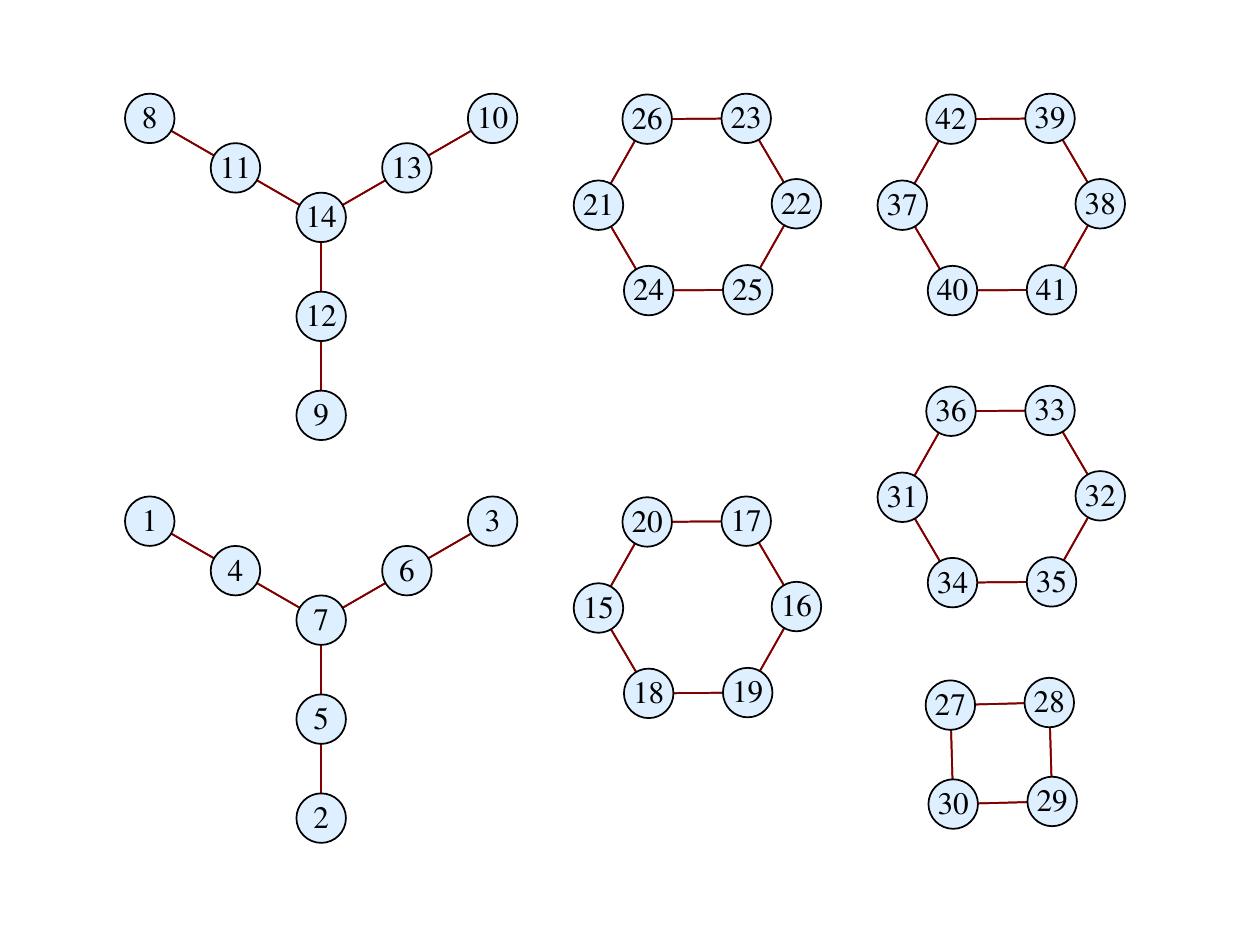} 
  \end{center}
\caption{Fundamental fusion graph of the quantum double of the binary tetrahedral}
 \label{fig:E6}
\end{figure}

The fusion graph of the fundamental irrep (numbered $4$), as displayed on figure \ref{fig:E6}, contains $7$ connected components: two copies of the affine $\widehat E_6$ Dynkin diagram, as expected from the McKay correspondence, four hexagons, and one square.

\noindent
The blocks of the symmetric and unitary matrix $24 \times S(I,J)$, for $I=1\ldots 7$, $I \leq J \leq 7$ are given below (we  set $\varphi = 1 + i \sqrt 3$ and $\psi = 1 - i \sqrt 3$): 
\smallskip

{\tiny 
$I=1$

$$
\left(
\begin{array}{ccccccc}
 1 & 1 & 1 & 2 & 2 & 2 & 3 \\
 1 & 1 & 1 & 2 & 2 & 2 & 3 \\
 1 & 1 & 1 & 2 & 2 & 2 & 3 \\
 2 & 2 & 2 & 4 & 4 & 4 & 6 \\
 2 & 2 & 2 & 4 & 4 & 4 & 6 \\
 2 & 2 & 2 & 4 & 4 & 4 & 6 \\
 3 & 3 & 3 & 6 & 6 & 6 & 9 \\
\end{array}
\right),\left(
\begin{array}{ccccccc}
 1 & 1 & 1 & 2 & 2 & 2 & 3 \\
 1 & 1 & 1 & 2 & 2 & 2 & 3 \\
 1 & 1 & 1 & 2 & 2 & 2 & 3 \\
 -2 & -2 & -2 & -4 & -4 & -4 & -6 \\
 -2 & -2 & -2 & -4 & -4 & -4 & -6 \\
 -2 & -2 & -2 & -4 & -4 & -4 & -6 \\
 3 & 3 & 3 & 6 & 6 & 6 & 9 \\
\end{array}
\right),\left(
\begin{array}{cccccc}
 4 & 4 & 4 & 4 & 4 & 4 \\
 -2 \psi  & -2 \psi  & -2 \psi  & -2 \psi  & -2 \psi  & -2 \psi  \\
 -2 \varphi  & -2 \varphi  & -2 \varphi  & -2 \varphi  & -2 \varphi  & -2 \varphi  \\
 -4 & -4 & -4 & -4 & -4 & -4 \\
 2 \psi  & 2 \psi  & 2 \psi  & 2 \psi  & 2 \psi  & 2 \psi  \\
 2 \varphi  & 2 \varphi  & 2 \varphi  & 2 \varphi  & 2 \varphi  & 2 \varphi  \\
 0 & 0 & 0 & 0 & 0 & 0 \\
\end{array}
\right)
$$

$$
\left(
\begin{array}{cccccc}
 4 & 4 & 4 & 4 & 4 & 4 \\
 -2 \varphi  & -2 \varphi  & -2 \varphi  & -2 \varphi  & -2 \varphi  & -2 \varphi  \\
 -2 \psi  & -2 \psi  & -2 \psi  & -2 \psi  & -2 \psi  & -2 \psi  \\
 -4 & -4 & -4 & -4 & -4 & -4 \\
 2 \varphi  & 2 \varphi  & 2 \varphi  & 2 \varphi  & 2 \varphi  & 2 \varphi  \\
 2 \psi  & 2 \psi  & 2 \psi  & 2 \psi  & 2 \psi  & 2 \psi  \\
 0 & 0 & 0 & 0 & 0 & 0 \\
\end{array}
\right),\left(
\begin{array}{cccc}
 6 & 6 & 6 & 6 \\
 6 & 6 & 6 & 6 \\
 6 & 6 & 6 & 6 \\
 0 & 0 & 0 & 0 \\
 0 & 0 & 0 & 0 \\
 0 & 0 & 0 & 0 \\
 -6 & -6 & -6 & -6 \\
\end{array}
\right),\left(
\begin{array}{cccccc}
 4 & 4 & 4 & 4 & 4 & 4 \\
 -2 \varphi  & -2 \varphi  & -2 \varphi  & -2 \varphi  & -2 \varphi  & -2 \varphi  \\
 -2 \psi  & -2 \psi  & -2 \psi  & -2 \psi  & -2 \psi  & -2 \psi  \\
 4 & 4 & 4 & 4 & 4 & 4 \\
 -2 \varphi  & -2 \varphi  & -2 \varphi  & -2 \varphi  & -2 \varphi  & -2 \varphi  \\
 -2 \psi  & -2 \psi  & -2 \psi  & -2 \psi  & -2 \psi  & -2 \psi  \\
 0 & 0 & 0 & 0 & 0 & 0 \\
\end{array}
\right)$$

$$\left(
\begin{array}{cccccc}
 4 & 4 & 4 & 4 & 4 & 4 \\
 -2 \psi  & -2 \psi  & -2 \psi  & -2 \psi  & -2 \psi  & -2 \psi  \\
 -2 \varphi  & -2 \varphi  & -2 \varphi  & -2 \varphi  & -2 \varphi  & -2 \varphi  \\
 4 & 4 & 4 & 4 & 4 & 4 \\
 -2 \psi  & -2 \psi  & -2 \psi  & -2 \psi  & -2 \psi  & -2 \psi  \\
 -2 \varphi  & -2 \varphi  & -2 \varphi  & -2 \varphi  & -2 \varphi  & -2 \varphi  \\
 0 & 0 & 0 & 0 & 0 & 0 \\
\end{array}
\right)
$$
}

{\tiny 
$I=2$

$$
\left(
\begin{array}{ccccccc}
 1 & 1 & 1 & -2 & -2 & -2 & 3 \\
 1 & 1 & 1 & -2 & -2 & -2 & 3 \\
 1 & 1 & 1 & -2 & -2 & -2 & 3 \\
 -2 & -2 & -2 & 4 & 4 & 4 & -6 \\
 -2 & -2 & -2 & 4 & 4 & 4 & -6 \\
 -2 & -2 & -2 & 4 & 4 & 4 & -6 \\
 3 & 3 & 3 & -6 & -6 & -6 & 9 \\
\end{array}
\right),\left(
\begin{array}{cccccc}
 4 & -4 & 4 & -4 & 4 & -4 \\
 -2 \psi  & 2 \psi  & -2 \psi  & 2 \psi  & -2 \psi  & 2 \psi  \\
 -2 \varphi  & 2 \varphi  & -2 \varphi  & 2 \varphi  & -2 \varphi  & 2 \varphi  \\
 -4 & 4 & -4 & 4 & -4 & 4 \\
 2 \psi  & -2 \psi  & 2 \psi  & -2 \psi  & 2 \psi  & -2 \psi  \\
 2 \varphi  & -2 \varphi  & 2 \varphi  & -2 \varphi  & 2 \varphi  & -2 \varphi  \\
 0 & 0 & 0 & 0 & 0 & 0 \\
\end{array}
\right),\left(
\begin{array}{cccccc}
 4 & -4 & 4 & -4 & 4 & -4 \\
 -2 \varphi  & 2 \varphi  & -2 \varphi  & 2 \varphi  & -2 \varphi  & 2 \varphi  \\
 -2 \psi  & 2 \psi  & -2 \psi  & 2 \psi  & -2 \psi  & 2 \psi  \\
 -4 & 4 & -4 & 4 & -4 & 4 \\
 2 \varphi  & -2 \varphi  & 2 \varphi  & -2 \varphi  & 2 \varphi  & -2 \varphi  \\
 2 \psi  & -2 \psi  & 2 \psi  & -2 \psi  & 2 \psi  & -2 \psi  \\
 0 & 0 & 0 & 0 & 0 & 0 \\
\end{array}
\right)
$$

$$
\left(
\begin{array}{cccc}
 6 & -6 & 6 & -6 \\
 6 & -6 & 6 & -6 \\
 6 & -6 & 6 & -6 \\
 0 & 0 & 0 & 0 \\
 0 & 0 & 0 & 0 \\
 0 & 0 & 0 & 0 \\
 -6 & 6 & -6 & 6 \\
\end{array}
\right),\left(
\begin{array}{cccccc}
 4 & -4 & 4 & -4 & 4 & -4 \\
 -2 \varphi  & 2 \varphi  & -2 \varphi  & 2 \varphi  & -2 \varphi  & 2 \varphi  \\
 -2 \psi  & 2 \psi  & -2 \psi  & 2 \psi  & -2 \psi  & 2 \psi  \\
 4 & -4 & 4 & -4 & 4 & -4 \\
 -2 \varphi  & 2 \varphi  & -2 \varphi  & 2 \varphi  & -2 \varphi  & 2 \varphi  \\
 -2 \psi  & 2 \psi  & -2 \psi  & 2 \psi  & -2 \psi  & 2 \psi  \\
 0 & 0 & 0 & 0 & 0 & 0 \\
\end{array}
\right),\left(
\begin{array}{cccccc}
 4 & -4 & 4 & -4 & 4 & -4 \\
 -2 \psi  & 2 \psi  & -2 \psi  & 2 \psi  & -2 \psi  & 2 \psi  \\
 -2 \varphi  & 2 \varphi  & -2 \varphi  & 2 \varphi  & -2 \varphi  & 2 \varphi  \\
 4 & -4 & 4 & -4 & 4 & -4 \\
 -2 \psi  & 2 \psi  & -2 \psi  & 2 \psi  & -2 \psi  & 2 \psi  \\
 -2 \varphi  & 2 \varphi  & -2 \varphi  & 2 \varphi  & -2 \varphi  & 2 \varphi  \\
 0 & 0 & 0 & 0 & 0 & 0 \\
\end{array}
\right)
$$
}

{\tiny 
$I=3$

$$
\left(
\begin{array}{cccccc}
 4 & 4 & -2 \varphi  & -2 \varphi  & -2 \psi  & -2 \psi  \\
 4 & 4 & -2 \varphi  & -2 \varphi  & -2 \psi  & -2 \psi  \\
 -2 \varphi  & -2 \varphi  & -2 \psi  & -2 \psi  & 4 & 4 \\
 -2 \varphi  & -2 \varphi  & -2 \psi  & -2 \psi  & 4 & 4 \\
 -2 \psi  & -2 \psi  & 4 & 4 & -2 \varphi  & -2 \varphi  \\
 -2 \psi  & -2 \psi  & 4 & 4 & -2 \varphi  & -2 \varphi  \\
\end{array}
\right),\left(
\begin{array}{cccccc}
 4 & 4 & -2 \varphi  & -2 \varphi  & -2 \psi  & -2 \psi  \\
 4 & 4 & -2 \varphi  & -2 \varphi  & -2 \psi  & -2 \psi  \\
 -2 \psi  & -2 \psi  & 4 & 4 & -2 \varphi  & -2 \varphi  \\
 -2 \psi  & -2 \psi  & 4 & 4 & -2 \varphi  & -2 \varphi  \\
 -2 \varphi  & -2 \varphi  & -2 \psi  & -2 \psi  & 4 & 4 \\
 -2 \varphi  & -2 \varphi  & -2 \psi  & -2 \psi  & 4 & 4 \\
\end{array}
\right),\left(
\begin{array}{cccc}
 0 & 0 & 0 & 0 \\
 0 & 0 & 0 & 0 \\
 0 & 0 & 0 & 0 \\
 0 & 0 & 0 & 0 \\
 0 & 0 & 0 & 0 \\
 0 & 0 & 0 & 0 \\
\end{array}
\right)
$$

$$
\left(
\begin{array}{cccccc}
 4 & 4 & -2 \varphi  & -2 \varphi  & -2 \psi  & -2 \psi  \\
 -4 & -4 & 2 \varphi  & 2 \varphi  & 2 \psi  & 2 \psi  \\
 -2 \psi  & -2 \psi  & 4 & 4 & -2 \varphi  & -2 \varphi  \\
 2 \psi  & 2 \psi  & -4 & -4 & 2 \varphi  & 2 \varphi  \\
 -2 \varphi  & -2 \varphi  & -2 \psi  & -2 \psi  & 4 & 4 \\
 2 \varphi  & 2 \varphi  & 2 \psi  & 2 \psi  & -4 & -4 \\
\end{array}
\right),\left(
\begin{array}{cccccc}
 4 & 4 & -2 \varphi  & -2 \varphi  & -2 \psi  & -2 \psi  \\
 -4 & -4 & 2 \varphi  & 2 \varphi  & 2 \psi  & 2 \psi  \\
 -2 \varphi  & -2 \varphi  & -2 \psi  & -2 \psi  & 4 & 4 \\
 2 \varphi  & 2 \varphi  & 2 \psi  & 2 \psi  & -4 & -4 \\
 -2 \psi  & -2 \psi  & 4 & 4 & -2 \varphi  & -2 \varphi  \\
 2 \psi  & 2 \psi  & -4 & -4 & 2 \varphi  & 2 \varphi  \\
\end{array}
\right)
$$
}

{\tiny
$ I = 4$

$$
\left(
\begin{array}{cccccc}
 4 & 4 & -2 \psi  & -2 \psi  & -2 \varphi  & -2 \varphi  \\
 4 & 4 & -2 \psi  & -2 \psi  & -2 \varphi  & -2 \varphi  \\
 -2 \psi  & -2 \psi  & -2 \varphi  & -2 \varphi  & 4 & 4 \\
 -2 \psi  & -2 \psi  & -2 \varphi  & -2 \varphi  & 4 & 4 \\
 -2 \varphi  & -2 \varphi  & 4 & 4 & -2 \psi  & -2 \psi  \\
 -2 \varphi  & -2 \varphi  & 4 & 4 & -2 \psi  & -2 \psi  \\
\end{array}
\right),\left(
\begin{array}{cccc}
 0 & 0 & 0 & 0 \\
 0 & 0 & 0 & 0 \\
 0 & 0 & 0 & 0 \\
 0 & 0 & 0 & 0 \\
 0 & 0 & 0 & 0 \\
 0 & 0 & 0 & 0 \\
\end{array}
\right),\left(
\begin{array}{cccccc}
 4 & 4 & -2 \psi  & -2 \psi  & -2 \varphi  & -2 \varphi  \\
 -4 & -4 & 2 \psi  & 2 \psi  & 2 \varphi  & 2 \varphi  \\
 -2 \psi  & -2 \psi  & -2 \varphi  & -2 \varphi  & 4 & 4 \\
 2 \psi  & 2 \psi  & 2 \varphi  & 2 \varphi  & -4 & -4 \\
 -2 \varphi  & -2 \varphi  & 4 & 4 & -2 \psi  & -2 \psi  \\
 2 \varphi  & 2 \varphi  & -4 & -4 & 2 \psi  & 2 \psi  \\
\end{array}
\right)
$$

$$
\left(
\begin{array}{cccccc}
 4 & 4 & -2 \psi  & -2 \psi  & -2 \varphi  & -2 \varphi  \\
 -4 & -4 & 2 \psi  & 2 \psi  & 2 \varphi  & 2 \varphi  \\
 -2 \varphi  & -2 \varphi  & 4 & 4 & -2 \psi  & -2 \psi  \\
 2 \varphi  & 2 \varphi  & -4 & -4 & 2 \psi  & 2 \psi  \\
 -2 \psi  & -2 \psi  & -2 \varphi  & -2 \varphi  & 4 & 4 \\
 2 \psi  & 2 \psi  & 2 \varphi  & 2 \varphi  & -4 & -4 \\
\end{array}
\right)
$$

{\tiny
$ I = 5$

$$
\left(
\begin{array}{cccc}
 12 & 0 & -12 & 0 \\
 0 & -12 & 0 & 12 \\
 -12 & 0 & 12 & 0 \\
 0 & 12 & 0 & -12 \\
\end{array}
\right),\left(
\begin{array}{cccccc}
 0 & 0 & 0 & 0 & 0 & 0 \\
 0 & 0 & 0 & 0 & 0 & 0 \\
 0 & 0 & 0 & 0 & 0 & 0 \\
 0 & 0 & 0 & 0 & 0 & 0 \\
\end{array}
\right),\left(
\begin{array}{cccccc}
 0 & 0 & 0 & 0 & 0 & 0 \\
 0 & 0 & 0 & 0 & 0 & 0 \\
 0 & 0 & 0 & 0 & 0 & 0 \\
 0 & 0 & 0 & 0 & 0 & 0 \\
\end{array}
\right)
$$
}

{\tiny
$I=6$

$$
\left(
\begin{array}{cccccc}
 4 & -4 & -2 \psi  & 2 \psi  & -2 \varphi  & 2 \varphi  \\
 -4 & 4 & 2 \psi  & -2 \psi  & 2 \varphi  & -2 \varphi  \\
 -2 \psi  & 2 \psi  & -2 \varphi  & 2 \varphi  & 4 & -4 \\
 2 \psi  & -2 \psi  & 2 \varphi  & -2 \varphi  & -4 & 4 \\
 -2 \varphi  & 2 \varphi  & 4 & -4 & -2 \psi  & 2 \psi  \\
 2 \varphi  & -2 \varphi  & -4 & 4 & 2 \psi  & -2 \psi  \\
\end{array}
\right),\left(
\begin{array}{cccccc}
 4 & -4 & -2 \psi  & 2 \psi  & -2 \varphi  & 2 \varphi  \\
 -4 & 4 & 2 \psi  & -2 \psi  & 2 \varphi  & -2 \varphi  \\
 -2 \varphi  & 2 \varphi  & 4 & -4 & -2 \psi  & 2 \psi  \\
 2 \varphi  & -2 \varphi  & -4 & 4 & 2 \psi  & -2 \psi  \\
 -2 \psi  & 2 \psi  & -2 \varphi  & 2 \varphi  & 4 & -4 \\
 2 \psi  & -2 \psi  & 2 \varphi  & -2 \varphi  & -4 & 4 \\
\end{array}
\right)
$$
}

{\tiny
$ I = 7$

$$
\left(
\begin{array}{cccccc}
 4 & -4 & -2 \varphi  & 2 \varphi  & -2 \psi  & 2 \psi  \\
 -4 & 4 & 2 \varphi  & -2 \varphi  & 2 \psi  & -2 \psi  \\
 -2 \varphi  & 2 \varphi  & -2 \psi  & 2 \psi  & 4 & -4 \\
 2 \varphi  & -2 \varphi  & 2 \psi  & -2 \psi  & -4 & 4 \\
 -2 \psi  & 2 \psi  & 4 & -4 & -2 \varphi  & 2 \varphi  \\
 2 \psi  & -2 \psi  & -4 & 4 & 2 \varphi  & -2 \varphi  \\
\end{array}
\right)
$$

\vfill 
\eject

\end{document}